\documentclass[12pt]{article}  
\usepackage{pudlatex,fullpage}
\usepackage{amsmath}

\newcommand\Rfn{\mbox{\it Rfn}}

\newcommand\Prf{\mbox{\it Prf}}
\newcommand\pf{\mbox{\it prf}}
\newcommand\rfn{\mbox{\it rfn}}

\title{Reflection principles, propositional proof systems, and theories
\\~\\
{\Large \it In memory of Gaisi Takeuti}}

\author{Pavel Pudl\'ak
\thanks{The author is supported by the project EPAC, funded by the Grant Agency of the Czech Republic under the grant agreement no. 19-27871X,  and the institute grant RVO: 67985840.}}

\begin{document}

\maketitle

\begin{abstract}
The reflection principle is the statement that if a sentence is provable then it is true. Reflection principles have been studied for first-order theories, but they also play an important role in propositional proof complexity. In this paper we will revisit some results about the reflection principles for propositional proofs systems using a finer scale of reflection principles. We will use the result that proving lower bounds on Resolution proofs is hard in Resolution. This appeared first in the recent article of Atserias and M\"uller~\cite{atserias-muller} as a key lemma and was generalized and simplified in some spin-off papers~\cite{garlik19,goos-et-al,garlik20}. We will also survey some results about arithmetical theories and proof systems associated with them. We will show a connection between a conjecture about proof complexity of finite consistency statements and a statement about proof systems associated with a theory.
\end{abstract}

\section{Introduction}

This paper is essentially a survey of some well-known results in proof complexity supplemented with some observations. In most cases we will also sketch or give an idea of the proofs. Our aim is to focus on some interesting results rather than giving a complete account of known results. We presuppose knowledge of basic concepts and theorems in proof complexity. An excellent source is Kraj\'{\i}\v{c}ek's last book~\cite{krajicek-proof}, where you can find the necessary definitions and read more about the results mentioned here. 

Proof complexity, as we view it, is not only the study of the complexity of propositional proofs, but also the study of first-order theories of arithmetic in connection with questions in computational complexity. In the seminal paper~\cite{cook} Stephen A. Cook defined a theory $PV$ and showed a close connection to the proof system \emph{Extended Resolution,} which is in a well-defined sense equivalent to the more familiar proof system  \emph{Extended Frege,} ($EF$). Similar connections have been shown between fragments of bounded arithmetic $T^i_2$, introduced in~\cite{buss}, and fragments of the quantified propositional calculus. In this paper we will continue this line of research in connection with the recent progress concerning the computational complexity of proof search in Resolution and some other weak proof systems.

The basic question about proof search is whether short proofs in a propositional proof system $P$ can be efficiently constructed if we know that they exist. More precisely this is the following problem: given a proposition $\phi$ such that there exists a proof of length $m$ in a proof system $P$, can we find a $P$-proof of $\phi$ in polynomial time? If this is possible, we say that $P$ is \emph{automatable}. A weaker property is \emph{weak automatability} where we only want to find a proof in a possibly stronger proof system $Q$. If $P$ is automatable, then it is possible, in particular, to \emph{decide} whether $\phi$ has a $P$-proof of length $\leq p(m)$, or does not have a proof of length $\leq m$ for some polynomial $p$.

One can show that all sufficiently strong proof systems are not weakly automatable iff there exists a disjoint {\bf NP} pair that is not separable by a set in {\bf P}. The latter condition follows, for instance, from the hypothesis that {\bf P}$\neq${\bf NP}$\cap${\bf coNP}. Therefore we believe that propositional proof systems are not-automatable except for some weak systems. The hardness assumptions about specific functions used in cryptography, such as the discrete logarithm, imply that already depth-$d$ Frege systems, for some small $d$, are not weakly automatable~\cite{bonet-et-al}. 

Our original motivation for this paper was the recent result of A. Atserias and M. M\"uller~\cite{atserias-muller} that Resolution is not automatable unless {\bf P}={\bf NP}. The essence of the proof is to define a polynomial reduction from the {\bf NP}-complete problem SAT to the problem to decide if a CNF formula has a short Resolution refutation. They define the reduction by mapping a CNF formula $\phi$ to a CNF formula $\rho_{\phi,m}$ expressing that $\phi$ has a Resolution refutation of length $m$ where $m=p(|\phi|)$ for a suitable polynomial $p$. Then they prove:
\ben
\item If $\phi$ is satisfiable, then $\rho_{\phi,m}$ has a Resolution refutation of length $q(m)$ where $q$ is some polynomial;
\item if $\phi$ is not satisfiable, then $\rho_{\phi,m}$ does not have a Resolution refutation of polynomial length.
\een 
Hence if there were a polynomial time algorithm for constructing Resolution proofs of length at most polynomially longer than the minimal ones, then it could be used to decide the satisfiability of $\phi$.%
\footnote{Atserias and M\"uller prove in fact a stronger lower bound in 2., namely there is no subexponential proof of $\rho_{\phi,n}$.} 

Our aim here is not to deal with automatizability, but rather address the natural question: for which propositional proof systems the conditions 1. and 2. are satisfied? The first condition is connected with the reflection principle. The reflection principle is the implication: if $\phi$ is provable, then $\phi$ is true. For every propositional proof system $P$, this can be formalized in the propositional calculus by a series of tautologies. It is well-known that some natural proof system, e.g., $EF$, prove these tautologies by polynomial length proofs. It is not difficult to see that if this is the case, then also formulas of the type $\rho_{\phi,n}$ have polynomial length proofs in the same proof system. Concerning condition 2., this is a question whether a proof system $P$ is able to prove superpolynomial lower bounds on $P$-proofs, i.e., its own proofs. The question is widely open, except for Resolution, and we consider it a more important than automatability.

This paper is organized as follows. 
After two preliminary sections, in Section 3, we will deal with reflection principles. So far only the full reflection principle and the consistency were studied in the context of the propositional calculus. We will introduce an intermediate concept of the \emph{local reflection principle}, which has only been studied in first order theories. Furthermore, this principle splits into two: the local reflection principle for tautologies and the local reflection principle for non-tautologies. With this finer scale in hand, we will revisit some previous results. Using the lower bound on Resolution proofs on unprovability in Resolution from~\cite{atserias-muller} we will show that Resolution proves efficiently the local reflection principle for non-tautologies, but not for tautologies. 

We know that many concrete proof systems prove their reflection principle, or at least the local reflection principle for non-tautologies, with polynomial length proofs. So in order to prove hardness of proof search in these systems, it would now suffice to prove lower bounds on proving lower bounds in these proof systems. But this is a nontrivial problem even for fairly weak proof systems and for strong systems one can only hope for some conditional results. We only know that a sufficiently strong proof system is not able to prove superlinear lower bounds on the proofs in a \emph{stronger} proof system.

At the end of Section 3 we will survey results on reflection principles for some concrete proof systems and state open problems.

In Section 4 we will study connections of first order theories and propositional proof system. For a theory $T$ extending some fragments of arithmetic, say $S^1_2$, one can define two associated proof systems. The first one, which may not always exist and which we call the weak proof system of $T$, is the strongest proof system $P$ whose soundness is provable in $T$. The second one, which always exists and we call it the strong proof system of $T$, is the proof systems in which the proofs are the first order proofs in $T$ of the statements that formulas are tautologies. The basic question is which theories have weak proof systems and which proof systems are weak, or strong,  proof systems of theories. 
One can show that $S^1_2+A$, where $A$ is a $\forall\Sigma^b_0$ sentence, always has the weak proof system, and a sufficiently strong proof system is the weak proof system of some theory iff it proves efficiently its reflection principle (more precisely, if $S^1_2$ proves this fact). Using a result about provability of finite consistency statements, Kraj\'i\v cek proved that the strong proof system of $T$ is equivalent to the weak proof system of $T+Con_T$, where $Con_T$ is the statement expressing the consistency of $T$; see~\cite{krajicek-proof}. We will connect the problem whether or not the weak proof system polynomially simulates the strong one with a conjecture about finite consistency statements.

\paragraph{Acknowledgment} I would like to thank Emil Je\v r\' abek, Erfan Khaniki, J\' an Pich, and especially Jan Kraj\' i\v cek for their remarks on the draft.


\section{Basic concepts}

\subsection{Propositional proof systems}

In this paper \emph{''a proof system''} will always mean \emph{a proof system for classical propositional logic}.
We will consider propositional proof systems in the sense of Cook and Rackhow~\cite{cook-reckhow}. According to this definition, a proof system is any polynomial time computable function $P$ that maps the set of strings $\{0,1\}^*$ onto the set of propositional tautologies {\it Taut}. The set {\it Taut}, of course, depends on the chosen basis of connectives and often we will only consider DNF tautologies. We will also study what is provable about proof systems in some first-order theories. Then it is important how the proof system is represented. Proof systems will always be represented by polynomial time algorithms, which we can formalize, say, by Turing machines, but we will also need that the theory in question can ``recognize'' some properties of the proof system. In order to prove nontrivial facts about proof systems, we will assume some conditions that guarantee that the proof systems are well-behaved. Thus we will also need that these conditions be provable in the theories in which we formalize our arguments about proof systems.  

We will prove some general statements about strong proof systems, at least as strong as Extended Frege systems, and survey some results about weak systems, specifically Resolution, Cutting Planes, and bounded depth Frege systems. We assume that the reader knows these basic systems. Unless stated otherwise, we will treat Resolution as a proof system for proving DNF tautologies, rather than refutation system for refuting unsatisfiable CNFs. 

Concerning the strong systems, it is more convenient to use 
Circuit Frege systems instead of the more familiar Extended Frege systems. Circuit Ferege systems are polynomislly equivalent to Extended Frege systems, Substitution Frege systems, Extended Resolution, and the sequent calculus with the extension rule {$ePK$}. The advantage of Circuit Frege systems is that one can use formalizations of polynomial time relations by circuits. If we have to formalize relations by formulas, we need additional \emph{extension variables} and \emph{extension axioms}. This is only a minor technicality, but it complicates notation, therefore we will use Circuit Frege proof system in this article.

\emph{Circuit Frege systems,} abbreviated $CF$, as defined by Je\v{r}\'abek~\cite{jerabek}, use circuits instead of formulas and the same rules as a Frege system with and additional rule that enables one to derive from a circuit $C$ a circuit $C'$  that unfolds in the same formula. Thus circuits are viewed as compressed forms of formulas that enable us to represent them more succinctly. Alternatively, we can compress circuits to a canonical incompressible form and identify circuits that have the same canonical form. 

It is not very important which representation one uses for representing formally circuits; we will assume that they are represented by straight-line programs and coded by 0--1 sequences. We will use the \emph{length of the sequence coding the circuit} as the measure of the complexity of circuits.

We say that \emph{a proof system $P$ extends a Circuit Frege system,} if $P$ is a a Circuit Frege system augmented with a polynomial time decidable set of sound axiom schemas $A$. Elements of $A$ are formulas or circuits that are tautologies and being schemas means that a proof may use \emph{any substitution instance} of a proposition $\alpha$ from $A$. When we formalize statements about a proof system $P$ that extend a Circuit Frege system in a theory $T$, we will always assume that $P$ is represented in such a way that $T$ proves that $P$ extends $CF$, i.e., $T$ proves that $P$-proofs are sequences of circuits etc.

It is well-known that $CF$ efficiently simulates the substitution rule~\cite{jerabek} and this also holds for extensions of $CF$ by schemas.
We will state a special case of this fact for further reference.

\bll{l-m}
Let $P$ be a proof system that extends $CF$. Suppose that $D$ is a $P$-proof of $\alpha(x_1\dts x_n)$ and let $\gamma_1\dts\gamma_n$ be circuits. Then one can construct a $P$-proof of $\alpha(\gamma_1\dts \gamma_n)$ in time polynomial in $|D|+\sum_i|\gamma_i|$. Furthermore, this is provable in $S^1_2$ provided that $P$ is formalized so that $S^1_2$ recognizes $P$ as an extension of $CF$.
\el
\bprf
Given a $P$-proof of $\alpha(x_1\dts x_n)$, we substitute  $\gamma_1\dts\gamma_n$ for $x_1\dts x_n$ in all circuits in $P$.
\eprf

In the context of proof systems extending $CF$, the word \emph{`proposition'} will mean a Boolean circuit (a special case of which is a formula). We will use it in particular when the circuit is supposed to express some truth. 

\subsection{Efficiently provable tautologies}

In propositional proof complexity theory we are mainly interested in asymptotical bounds. The typical question is whether or not a given sequence of tautologies has polynomial size proofs in a given proof system. We are often also interested whether the proofs can be efficeintly constructed and verified in a weak theory. 
Thus given a proof system $P$ and a sequence of tautologies $\{\phi_n\}$ we may consider three basic properties of increasing strength:
\ben
\item $\{\phi_n\}$ have $P$-proofs of polynomial length,
\item $P$-proofs of $\{\phi_n\}$ can be constructed in polynomial time,
\item $S^1_2$ proves that $\{\phi_n\}$ are provable in $P$.
\een
The third condition implies the first two if $\{\phi_n\}$ is 
constructible in polynomial time,
which we will usually assume. In the propositions and definitions in this paper we will use only one of these properties, mostly the third one, which does not mean that they do not have versions with the other two properties. We leave to the reader to prove the other versions when it is possible.
We will abbreviate the last two properties by 
\ben
\item[2.] $P$ $p$-proves $\{\phi_n\}$, and
\item[3.] $P$ provably $p$-proves $\{\phi_n\}$.
\een
If  $\{\phi_n\}$ is constructible in polynomial time, we will also say that 
\ben
\item[2.] $P$ $p$-proves  $\{\phi_n\}\to\{\psi_n\}$, 
\een
if $P$-proofs of $\{\psi_n\}$ can be constructed in polynomial time from  $\{\phi_n\}$.

\subsection{Theories}

Although some results could be generalized to a large class of theories, we will only consider finite extensions of $S^1_2$, the fragment of Bounded Arithmetic introduced by Buss~\cite{buss}. In $S^1_2$ polynomial time computations have natural formalization. Sometimes it convenient to have terms, rahter than formulas, for polynomial time algorithms. Then we will assume that the theories contain $S^1_2(PV)$ as defined in~\cite{buss}. In this theory every polynomial time function can be represented by a term.
We will assume that first-oder proofs are formalized by the standard Hilbert-style calculus. 

Bounded formulas are classified into classes $\Sigma^b_i$ and $\Pi^b_i$ according to the number of alternations of bounded existential and universal quantifier, ignoring sharply bounded quantifiers. A sharply bounded quantifier is a quantifier where the variable is bounded by a term of the form $|t|$, with $t$ an arbitrary arithmetical term and $|\dots|$ denoting the length, which is approximately the binary logarithm. Thus sharp bounds restrict the range of quantification to a polynomial size domain. In particular, a $\Sigma^b_1$ is a formula in prenex form with quantifiers of the form $\exists x\leq t$ and $\forall y\leq|s|$, where $t,s$ are terms (or a formula equivalent to it). We will denote by $\forall\Sigma^b_0$ formulas in the prenex form starting with one, or several universal quantifiers, folowed by sharply bounded quantifiers. A $\Pi_1$ formula is a prenex formula starting with unbounded universal quantifiers followe by bounded quantifiers.

We will use two basic results about arithmetical theories.

\bt[Parik's and Buss's Theorems, \cite{parikh,buss}]
~

\ben
\item Let $T=S^1_2+A$ where $A$ is a set of $\Pi_1$ sentences in the language of $S^1_2$, let $\phi(x,y)$ be a bounded formula with two variables, and suppose that $T\vdash\ \forall x\exists y.\phi(x,y)$. Then there exists a term $t(x)$ such that 
$T\vdash\ \forall x\exists y(y\leq t(x)\wedge\phi(x,y))$.%
\footnote{This is not the original form of Parik's Theorem and it is certainly not the most general form.}
\item Let $T=S^1_2+A$ where $A$ is a set of $\forall\Sigma^b_0$ sentences, let $\phi(x,y)$ be a $\Sigma^b_1$ formula with two variables, and suppose that $T\vdash\ \forall x\exists y.\phi(x,y)$. Then there exists a polynomial time computable function $f(x)$ such that 
$\N\models\ \forall x.\phi(x,f(x))$. Furthermore, if we extend $T$ with $PV$ to $T(PV)$ (i.e., $S^1_2(PV)+A$), then $T(PV)$ proves $\forall x.\phi(x,t(x))$ for some $PV$ term $t(x)$.
\een
\et

For a true $\forall\Sigma^b_0$ sentence $A$ of the form $\forall x.\alpha(x)$ with $\alpha\in\Sigma^b_0$, one can construct a sequence of polynomial length tautologies $[\![A]\!]_n$ that express that $\alpha(x)$ is satisfied for all $x$ of length $\leq n$. If we express $[\![A]\!]_n$ as a circuit, then it has $n$ variables; if it is a formula, then the number of variables bounded by a polynomial. 
The following is a fundamental result about theories and propositional proof system.

\bt[\cite{buss,cook}]\label{t-cook}
If $S^1_2$ proves a $\forall\Sigma^b_0$ sentence $A$, then $CF$ provably $p$-proves $[\![A]\!]_n$.
\et
This theorem has been extended to a number of theories and proof systems (see~\cite{cook-nguyen}), for example to fragments of Bounded Arithmetic and  fragments of the sequent calculus for quantified propositions~\cite{krajicek-pudlak-quantified}. We will say more about it in Section~4.

\subsection{Propositional proof systems vs. theories}

The last theorem shows that we can view provability of $\forall\Sigma^b_0$ sentences as a uniform way of proving tautologies. This is not the only possible way to represent uniform provability; e.g., the concept of $p$-provability can be viewed as being  intermediate between the nonuniform provability and uniform provability. The concept of provably $p$-provability is very close to the uniform provability as defined by provability of $\forall\Sigma^b_0$ sentences in first-order theories. In \cite{cook}, Cook proved the original version of Theorem~\ref{t-cook} for the equational theory $PV$, which is another intermediate step between propositional proof systems and first-order theories. Gaisi Takeuti introduced restricted types of first-order proofs that were close to propositional proofs. His aim was to use G\"odel's second incompleteness theorem to prove separations between fragments of Bounded Arithmetic. One of such concepts appeared in a joint work with Jan Kraj\'i\v cek~\cite{krajicek-takeuti}.


\section{Reflection principles and soundness}

\subsection{Definitions and general facts}

We are interested in reflection principles in propositional calculus, but we will start with reflection principles for first order theories as a paradigm. 
(For a survey of reflection pricniples in arithmetical theories, see Smorynski~\cite{smorynski-HB}.)
The \emph{local reflection principle} for a theory $T$ is the statement
\[
\mbox{\it for all $x$, if $x$ is a $T$-proof of $\phi$, then }\phi.
\]
The principle is called \emph{local} because it is stated for a single proposition. Usually we study a \emph{schema}, i.e., a set of such sentences for all propositions from some class of sentences $\cal C$. Given a class of sentences, we can also state the \emph{uniform reflection principle} for a theory $T$ and for class $\cal C$:
\[
\mbox{\it for all $x$ and $y$, if $y\in{\cal C}$ and $x$ is a $T$-proof of $\phi$, then  $y$ is true.}
\]
In order for this principle to be stated in the language of $T$, it must be possible to define the truth of sentences in $\cal C$. By Tarski's Theorem (which is easily provable using the fix-point lemma) it is not possible to define truth for all formulas. In particular, in an arithmetical theory we can define truth for classes $\Pi_n$, but not for all arithmetical formulas. We will see that this is different in the propositional calculus. What is also different is that in the propositional calculus we can only speak about proofs up to some length.

Now we want to state reflection principles for a propositional proof system $P$ in the propositional calculus. A natural way to do it is first to state it as an arithmetical formula of the form $\forall\Sigma^b_0$ and translate the formula into a sequence of propositions that express the principle for proofs up to length $n$. 

For a propositional proof system $P$ and a proposition $\phi$, we will denote by $\mbox{\it LRfn}_{P,\phi}$ a siutable natural formalization of the local reflection principle in $S^1_2$. In principle we can translate $\mbox{\it LRfn}_{P,\phi}$ to the propositional calculus as a sequence of propositions for every length $m$ of a proof, but for typical calculi the proofs are at most exponentially long in the length of the formula, so it does not make sense to state it for larger lengths. We will denote by $\mbox{\it lrfn}_{P,\phi,m}(\vec x)$ a suitable circuit expressing the principle for proofs $\vec x$ of length at most $m$, where $\vec x$ is a string of propositional variables $x_1\dts x_n$ representing a proof. 
The length of the circuit  $\mbox{\it lrfn}_{P,\phi,m}(\vec x)$ is bounded by a polynomial in $|\phi|$ and $m$. We denote by $|\phi|$ the length of the bit representation of $\phi$. If we cannot use circuits in a proof system in which we want to formalize the local reflexion principle, $\mbox{\it lrfn}_{P,\phi,m}$ will be a formula with additional extension variables.

Informally, the \emph{global reflection principle} for a proof system $P$ is the statement
\[
\mbox{\it for all $x$ and $y$, if $x$ is a $P$-proof of $y$, then $y$ is true}
\]
where \emph{`$y$ is true'} means that $y$ is a tautology. Again we can first state the principle in an arithmetical theory, namely $S^1_2$, and then translate it to the propositional calculus. 
The relation \emph{`$x$ is a $P$-proof of $y$'} is formalized by a $\Sigma^b_1$ formula $\Prf_P(x,y)$.  
The property that a number (or a bit string) encodes a tautology can be easily formalized by a strict $\Pi^b_1$ formula, which we will denote by $Taut(y)$.%
\footnote{A strict $\Pi^b_1$ formula is a formula that starts with one, or several, bounded universal quantifiers followed by a $\Sigma^b_0$ formula, i.e. formula with only sharply bounded quantifiers.}
Thus the global reflection principle for a proof system $P$ is formalized by the following $\forall\Sigma^b_0$ formula: 
\[
\Rfn_P:= \forall x,y(\mbox{\it Prf}_P(x,y)\to Taut(y)).
\]
In the rest of the paper we will omit the specification \emph{`global'} unless we need to stress the difference between global and local principles. It should be noted that the reflection principle for $P$ expresses the \emph{soundness} of~$P$; in this paper we will only refer to it by the term \emph{reflection principle.}

In the propositional calculus we can define circuits $\mbox{\it prf}_{m,n}(\vec x,\vec y)$ that express the relation \emph{`$x$ is a $P$-proof of $y$'} where $|x|=m$ and $|y|=n$,
and 
$sat_n(\vec y,\vec z)$ that express that the circuit encoded by $\vec y$ is satisfied by assignment $\vec z$ where $|y|=|z|=n$ ($z$ may have more bits than the number of variables of the proposition encoded by $x$). The propositional formalization of the reflection principle is the set of circuits
\[
\mbox{\it rfn}_{P,m,n}(\vec x,\vec y,\vec z):=\mbox{\it prf}_{P,m,n}(\vec x,\vec y)\to sat_n(\vec y,\vec z)
\]
for all $m,n\in\N$.
Since the reflection principle for $P$ is stated for all circuits, it is the same thing as the \emph{soundness} of $P$. 

Using the formalization of provability relation in $P$ we can state the local reflection principle $\mbox{\it lrfn}_{P,\phi,m}(\vec x)$ more explicitly by
\[
\mbox{\it prf}_{P,m,n}(\vec x,\lceil\phi\rceil)\to\phi,
\]
where $n$ is the length of $\phi$ and $\lceil\phi\rceil$ is the bit string representing $\phi$. 

If a proof $P$ system satisfies some basic properties, then the global reflection principle implies the local principle using short $P$-proofs.

\begin{fact}\label{f1}
Suppose that $P$ allows substitution of truth constants and $P$-proofs of propositions $\phi(\vec z)\equiv sat_n(\lceil\phi\rceil,\vec z)$ can be constructed in polynomial time. Then $P$-proofs of $\mbox{\it lrfn}_{P,\phi,m}(\vec x,\vec u)$ can be constructed in polynomial time from proofs of  $\mbox{\it rfn}_{P,m,n}(\vec x,\vec y,\vec z)$ where $m=|\phi|$  
(polynomial time is in $m$ and $n$).
\end{fact}

The condition about $\phi(\vec z)\equiv sat_n(\lceil\phi\rceil,\vec z)$ is not very restrictive provided that we formalize $sat_n$ in natural way. We can prove the equivalence by induction on the complexity of $\phi$. Then looking at this proof we can easily see that it actually gives us a polynomial time algorith to construct a $P$-proof of the equivalence.

The opposite implication in Fact~\ref{f1} is probably not true in general, but we do not have any example of a proof system with short proofs of all the instances of the local reflection principle and only long proofs of the global reflections principle. On the other hand, in sufficiently strong proof systems, in particular in proof systems that extend $CF$, already one particular instance of the local reflection principle implies the global reflection principle by short proofs. Namely, if we apply the local reflection principle with the formula $\bot$ (the constant representing contradiction), then it expresses the consistency of the proof system:
\[
\mbox{\it prf}_{P,m,n_0}(\vec x,\lceil\bot\rceil)\to\bot,
\]
which is equivalent to $\neg\mbox{\it prf}_{P,m,n_0}(\vec x,\lceil\bot\rceil)$, where $n_0$ is the length of the representation of $\bot$. We will denote this proposition by $con_{P,m}$.

\bpr 
Suppose that $S^1_2$ proves that $P$ is an extension of $CF$. Then $S^1_2$ proves: if $P$ is consistent, then it is sound.
\epr
\bprf
Suppose $D$ is a proof a formula $\alpha(\vec x)$. If $\alpha(\vec x)$ is not a tautology, then for some assignment $\vec a$, $\alpha(\vec a)$ is false. Since $P$ is an extension of $CF$, $P$ proves $\alpha(\vec a)$ (with the proof obtained from $D$ by substituting $\vec a)$ and it also proves $\neg\alpha(\vec a)$, which is a contradiction. So if $P$ is consistent, then $\alpha(\vec x)$ must be a tautology.

The above argument is elementary and the construction of the contradiction can be done in polynomial time. Hence the argument can be formalized in $S^1_2$. 
\eprf

\begin{corollary}\label{c3.2}
Let $P$ be an extension of $CF$. Then $P$-proofs of the propositions $\mbox{\it rfn}_{P,m,n}(\vec x,\vec y,\vec z)$ expressing the reflection principle for $P$  can be constructed in polynomial time from $P$-proofs of propositions $con_{P,m}$ expressing the consistency of $P$. 
\end{corollary}
Using our abbreviations, the conclusion can be stated as: 
\[
\mbox{$P$ $p$-proves $\{con_{P,m}\}\to\{\mbox{\it rfn}_{P,m,n}(\vec x,\vec y,\vec z)\}$.}
\]
\bprf
According to the previous proposition, $S^1_2$ proves
\[
\forall z\neg \Prf_P(z,\lceil\bot\rceil)\to\forall x,y(\Prf_P(x,y)\to \forall u.Sat(y,u)).
\]
By Buss's Theorem, there exists a $PV$ term $t$ such that 
\[
S^1_2(PV)\ \vdash\ \forall x,y,u(\neg \Prf_P(t(x,y,u),\lceil\bot\rceil)\to(\Prf_P(x,y)\to Sat(y,u)).
\]
Hence $CF$-proofs of propositional translations 
\bel{e-rf}
\neg \pf_{P,m_1,n_0}(\vec\gamma,\lceil\bot\rceil)\to(\pf_{P,m_2,n}(x,y)\to sat_n(y,u))
\ee
can be constructed in polynomial time, where $\vec\gamma$ is a string of circuits representing $t$. Since $t$ represents a polynomial time computable function, $|t(x,y,u)|$ is polynomial in $|x|,|y|,|u|$, in fact, polynomial in $|x|$ and $|y|$, because $|u|\leq|x|$. This means that $m_1$ is polynomial in $m_2$ and $n$. Hence we get $P$-proofs of (\ref{e-rf}) in polynomial time where the polynomial bound is in the length of $\rfn_{P,m_2,n}$. 
Lemma~{\ref{l-m}} guaranties that we can get $P$-proofs of $\neg\pf_{P,m_1,n_0}(\vec\gamma,\lceil\bot\rceil)$ from instances of $con_{P,n}$ in polynomial time. Finally, we get proofs of the instances the reflection principle for $P$ using modus ponens.
\eprf

We can furthermore distinguish two special cases of the local reflection principle: the \emph{local reflection principle for tautologies} and the \emph{local reflection principle for non-tautologies}. Consider the local reflection principle for a proposition $\phi$ stated as a disjunction
\[
\phi\vee\neg\mbox{\it prf}_{P,m,n}(\vec x,\lceil\phi\rceil).
\]
Suppose it has a short $P$-proof. If $\phi$ is a non-tautology, we can substitute a falsifying assignment into $\phi$ and thus get a short proof of the fact that $\phi$ does not have proof of length $\leq m$. Now suppose that $\phi$ is a tautology and, moreover, $P$ has the \emph{feasible disjunction property}. The latter property means that given a $P$-proof of a disjuntion of propositions with disjoint sets of variables, one can construct in polynomial time a $P$-proof of one of the disjuncts.
So if the reflection principle has a short $P$-proof, then either $\phi$ has a short $P$-proof, or the fact that it does not have proof of length $\leq m$ has a short proof. 

We now state these facts formally.

\bpr
Suppose a proof system $P$ has the property that from a $P$-proof of a disjucntion $\phi(\vec x)\vee\psi(\vec y)$, where $\phi(\vec x)$ is a nontautology, one can derive by a polynomial length $P$-proof of $\psi(\vec y)$, then TFAE
\ben
\item the instances of the local reflection principle for non-tautologies have polynomial length proofs, and
\item the propositions of the form $\neg\mbox{\it prf}_{P,m,n}(\vec x,\lceil\phi\rceil)$ for non-tautologies have polynomial length $P$-proofs.
\een
\epr

\bpr
Suppose a proof system $P$ has the feasible disjunction property and proves the instances of its local reflection principle for tautologies by proofs of polynomial lengths. Let  $\{\phi_n\}$ be a sequence of tautologies such that
$|\phi_n|\leq n$, and the shortest $P$-proof of $\phi_n$ has length $\geq m(n)$ where $m(n)\geq n$ is some function. Then $P$ can prove lower bounds on $P$ proofs of $\{\phi_n\}$ of the form $\geq m(n)$ using proofs whose length is polynomial in $m(n)$.
\epr
This fact was mentioned in~\cite{disjointNP} for the reflection principle. Here we use a weaker assumption, the local reflection principle for tautologies.

\subsection{Non-automatability of proof search and provable lower bounds}

In order to use the approach of Atserias and M\"uller for proving that proof search is hard for $P$,%
\footnote{This is certainly not the only possible way of proving hardness of proof search. The formulas used in the spin-off papers that proved hardness of proof search for  proof systems stronger than resolution still used  $\neg\mbox{\it prf}_{Res,p(n),n}(\vec x,\lceil\phi\rceil)$, it only was suitably lifted.} 
it suffices to prove for some polynomials $p(x)$ and $q(y)$ for every $n$, that
\ben
\item $\neg\mbox{\it prf}_{P,p(n),n}(\vec x,\lceil\phi\rceil)$ has $P$-proofs of length $q(n)$ when $\phi$ is \emph{not} a tautology and 
\item $\neg\mbox{\it prf}_{P,m,n}(\vec x,\lceil\phi\rceil)$  does not have polynomial length proofs if $\phi$ is a tautology.
\een

The first condition is satisfied when $P$ has polynomial length proofs of the instances of its local reflection principle for non-tautologies. In the next section we will show that there are many natural proof systems that prove their global reflection principle, hence this condition is satisfied for them. 
We will also mention lower bounds, i.e., condition 2, that have been proved for some weak systems.
Whether or not the second condition is satisfied by strong proof systems is not clear. We are only able to prove a weaker statement, which is not sufficient for non-automatibility. First we need a lemma.

\bll{l5.1}
Suppose a proof system $Q$ extending $CF$ is not sound. Then it proves every proposition using a linear length proof.
\el
\bprf
Let $Q$ prove $\alpha(\vec x)$ and suppose there is an assignment $\vec a$ that falsifies $\alpha(\vec x)$. Since $CF$ is complete, it proves $\neg\alpha(\vec a)$. Given a $Q$-proof of $\alpha(\vec x)$, substitute $\vec a$ for $\vec x$. Thus we obtain a $Q$-proof of $\alpha(\vec a)$. Hence $Q$ proves a contradiction, from which any proposition can be derived by a $CF$-proof of linear length.
\eprf

\bprl{p-unprov}
Let $P$ and $Q$ be extensions of $CF$ such that $P$ provably $p$-proves its reflection principle.
Let $\{\alpha_n\}$ be a polynomial time computable sequence of tautologies with $|\alpha_n|=n$. Then there exists a constant $c$ such that if $P$ provably $p$-proves a lower bound $>cn$ on $Q$-proofs of $\alpha_n$ for some $n_0$, then  $P$ polynomially simulates~$Q$.
\epr
Proving lower bound $cn$ on $Q$-proofs of $\alpha_n$  means provably $p$-proving $\neg\pf_{Q,cn,n}(x,\lceil\alpha_n\rceil)$.
The idea of the proof is simple: if $P$ provably $p$-proves a lower bound, then it also must prove the consistency of $Q$, hence also the reflection principle for $Q$.
We postpone the proof of this proposition until we elaborate connections between propositional proof systems and first order theories.

There are two weak points in Proposition~\ref{p-unprov}. First, it only refers to \emph{provable} $p$-provability, while we would like to have mere $p$-provability. Second, it only shows unprovablity for proof systems strictly stronger than $P$. Can any of the two weaknesses be removed?
 
The impossibility of proving superlinear lower bound resembles the situation with lower bounds on circuit complexity of concrete Boolean functions where we are only able to prove a lower bound $cn$ for some constant slightly larger than $3$. Is it only a coincidence?

\subsection{Reflection principles in concrete proof systems}

\subsubsection{Resolution and Cutting Planes}

We will denote the Resolution proof system by $Res$. Recall that in this article we view Resolution as proof system for proving DNF tautologies. We will use the following important result.

\bl[Atserias-M\"uller~\cite{atserias-muller}, Garl\'ik~\cite{garlik19}]\label{l39}
There exist constants $c>1$ and $\delta>0$ such that for every $n$, $m\geq n^c$, and proposition $\phi$, $|\phi|=n$, there is no Resolution proof 
\[
\neg\mbox{\it prf}_{Res,m,n}(\vec x,\lceil\phi\rceil)
\]
of length $n^\delta$.
\el
This is a generic statement of two different lemmas, because the encoding of Resolution proofs in the two cited papers are different. 

\begin{corollary}
There is no subexponential upper bound on the proofs of the $Res$-local reflection principle for tautologies in $Res$. 
Specifically, $Res$-proofs of the local reflection principle for $\alpha_n$, $|\alpha_n|=n$, have exponential length for {any sequence of tautologies $\{\alpha_n\}$ for which we have an exponential lower bound.}
\end{corollary}

\bprf
Suppose 
\[
\neg\mbox{\it prf}_{Res,m,n}(\vec x,\lceil\alpha_n\rceil)\vee \alpha_n,
\]
where $m\geq n^c$, has a proof of length $\leq S$. By the feasible disjunction property of $Res$, either $\neg\mbox{\it prf}_{Res,m,n}(\vec x,\lceil \alpha_n\rceil)$ has a proof of length $S^{O(1)}$, or $\alpha_n$  has a proof of length $S^{O(1)}$. Since both propositions have only exponentially large proofs, $S$ must be exponential.
\eprf

Since we do have sequences of tautologies that require exponential size proofs in Resolution, we also have such a lower bound on the \emph{local reflection principle for tautologies}. This implies an exponential lower bound for the \emph{global reflection principle} in resolution. 
On the other hand, one can construct polynomial length $Res$-proofs of $\neg\mbox{\it prf}_{Res,m,n}(\vec x,\lceil\phi\rceil)$ for non-tautologies $\phi$. This was proved in~\cite{disjointNP} and reproved in \cite{atserias-muller}. Hence the instances of the local reflection principle for non-tautologies have polynomial length proofs in $Res$. We state it for further reference.
\bpr\label{p39}
There exists a polynomial $p$ such that for every nontautology $\phi$, and $m\geq n$, $\neg\mbox{\it prf}_{Res,m,n}(\vec x,\lceil\phi\rceil)$ has a Resolution proof of length at most $p(m)$.
\epr
\bpri
In this proof we will treat Resolution as a refutation system. If $\phi$ is a non-tautology, then $\neg\phi$ is satisfiable. Suppose we have an assignment that satisfies all initial clauses in a refutation of $\phi$. We will gradually prove that every clause in the proof is satisfied, which produces a contradiction, because the last clause is empty. The whole point is that one can express the fact {\it ``a clause is satisfied by a given assignment''} by a clause; therefore this argument can be performed in Resolution.

In more detail, let the $j$th clause $C_j$ be represented by variables $y_{e,i,j}$ where the intended meaning is that $x_i$ ($\neg x_i$) is present in the clause if $e=1$ (respectively, $e=0$). Then the fact that the clause is satisfied by an assignment $(e_1\dts e_n)$ is expressed by 
\[
y_{e_1,1,j}\vee\dots\vee y_{e_n,n,j}.
\]
\epri

The lower bound on the reflection principle was proved already in 2004 by Atserias and Bonet~\cite{atserias-bonet}. It is based on an idea different from the one presented at the beginning of this section. We will explain that proof using a new concept.

\bdf
A disjunction $\phi(\vec p,\vec x)\vee\psi(\vec p,\vec y)$ with $\vec x$ and $\vec y$ disjoint strings of variables  is \emph{friendly for a proof system $P$} if
\ben
\item for every assignment $\vec p\vec x:=\vec a\vec b$ that falsifies $\phi$, $\psi(\vec a,\vec y)$ has a polynomial size $P$-proof, 
\item for every assignment $\vec p\vec y:=\vec a\vec c$ that falsifies $\psi$,  $\phi(\vec a,\vec x)$ has a polynomial size $P$-proof.
\een
\edf
More precisely, friendliness should be defined for a sequence of disjunctions 
$\{\phi_n(\vec p,\vec x)\vee\psi_n(\vec p,\vec y)\}$ and using $p$-provability in $P$ (we leave it to the reader).

\medskip
We say that $\phi(\vec p,\vec x)\vee\psi(\vec p,\vec y)$ is \emph{semi-friendly} for a proof system $P$ if condition 1. is satisfied, while condition 2. may fail. For proving a lower bound on the reflection principle, it would suffice to use semifriendly disjunctions, but 
we will use the stronger, more natural concept.

Here are some examples of friendly disjunctions.
\bpr
\ben
\item Every $p$-provable disjunction in any proof system.
\item $Clique_{k+1}$-$Coloring_k$ tautology is freindly in Cutting Planes.
\item Reflection principle is friendly in Resolution.
\een
\epr

\bpri
1. is trivial.

\medskip
\noindent 2. The tautology has the form
\[
\neg Clique_{k+1}(K,G)\vee \neg Coloring_k(\chi,G).
\]
where $K$ stands for a $k+1$-clique, $G$ for a graph, and $\chi$ for a $k$-coloring. One can show:
\bi
\item For a fixed graph $G$ and clique $K$, $\neg Coloring_k(\chi,G)$ follows from $PHP^{k+1}_k$.
\item For a fixed graph $G$ and coloring $\chi$, $\neg Clique_{k+1}(K,G)$ follows from $PHP^{k+1}_k$.
\ei

\medskip
\noindent 3. The tautology has the form
\[
\neg\mbox{\it prf}_{P,m,n}(\vec x,\vec y)\vee sat_n(\vec y,\vec z)
\]
One can show:
\bi
\item Given a $Res$-proof $\Pi$ of $\phi$, we also have a $Res$-proof of the equivalent formula $sat_n(\lceil\phi\rceil,\vec z)$, because $\phi\equiv sat_n(\lceil\phi\rceil,\vec z)$ is $p$-provable in $Res$.
\item Given $\phi$ and a falsifying assignment $\vec a$, we get a polynomial size proof of $\neg\mbox{\it prf}_{P,m,n}(\lceil\phi\rceil,\vec y)$ by Proposition~\ref{p39}.
\ei
\epri

We will first sketch the lower bound on the reflection principle in the Cutting Plane proof system. 

\bt[\cite{disjointNP}]
{$CP$-proofs of the reflection principle of $CP$ have size $\geq 2^{n^\epsilon}$ for some $\epsilon>0$.}
\et
\bpri
The idea is, roughly speaking, to construct a monotone polynomial reduction of the {\it Clique-Coloring} disjoint {\bf NP}-pair to the canonical pair of $CP$, which is, essentially, the pair defined by the reflection principle of $CP$. The reduction is defined by the mapping
\[
G\ \mapsto\ \neg Clique_{k+1}(\vec x,G).
\]
We observe that if $G$ has a $k$-coloring, then we can construct a proof of $\neg Clique_{k+1}(\vec x,G)$ using $PHP^{k+1}_k$. Here we use the fact that the Clique-Coloring disjunction is (semi)friendly in Cutting Planes. Let $p$ be a polynomial bound on such proofs.
Thus we map the disjoint {\bf NP} pair
\[
\left(\{G\ |\  Clique_{k+1}(\vec x,G) \}, \{G\ |\ Color_k(\vec y,G)\}\right)
\]
to the disjoint {\bf NP} pair
\[
\left(\{G\ |\  Clique_{k+1}(\vec x,G) \}, \{G\ |\ \neg Clique_{k+1}(\vec x,G) \mbox{ has a $CP$-proof of size }p(|G|)\}\right).
\]
From the fact that the first pair cannot be separated by subexponential monotone circuits, we get that the second one also cannot.
Using the monotone interpolation for $CP$ we get an exponential lower bound on $CP$ proofs of the disjunctions
\bel{e-negcl}
\neg Clique_{k+1}(\vec x,\vec y)\vee\neg\pf_{CP,p(n),n}(\vec z,\lceil\neg Clique_{k+1}(\vec x,\vec y)\rceil(\vec y))
\ee
that express the disjointness of the latter pair. 
The expression $\lceil\neg Clique_{k+1}(\vec x,\vec y)\rceil(\vec y)$ represents a circuit that given an assignment $\vec y:=\vec a$ of $0$s and $1$s, produces the string that is a code of the proposition $Clique_{k+1}(\vec x,\vec a)$.%
\footnote{In logic various different notations are used for the G\"odel numbers as functions of some variables. E.g., Smorynski~\cite{smorynski-HB} would write $\lceil\phi(x,\dot{y})\rceil$ for what we would denote by  $\lceil\phi(\vec x,\vec y)\rceil(\vec y)$.}
Hence the propositional variables of\\ $\neg\pf_{CP,p(n),n}(\vec z,\lceil\neg Clique_{k+1}(\vec x,\vec y)\rceil(\vec y))$ are only $\vec z$ and $\vec y$. Surely, in $CP$ we cannot represent nontrivial circuits, unless we use extension variables and extension axioms, but if we choose a suitable encoding, the circuit will be trivial---just constants and variables. Thus $\lceil\neg Clique_{k+1}(\vec x,\vec y)\rceil(\vec y)$ is simply a string of truth constants and variables $\vec y$.

Hence the disjunction~(\ref{e-negcl}) is, essentially, the reflection principle restricted to formulas of the form $\neg Clique_{k+1}(\vec x,\vec y)$. To get the special case of the reflection principle formally correct, we only need to replace $\neg Clique_{k+1}(\vec x,\vec y)$ with the polynomially equivalent formula $Sat_n(\lceil\neg Clique_{k+1}(\vec x,\vec y)\rceil,\vec u)$.
\epri

\bt[Atserias-Bonet \cite{atserias-bonet}]
$Res$-proofs of the reflection principle of $Res$ have size $\geq 2^{n^\epsilon}$ for some $\epsilon>0$.
\et
\bpri
The idea is 
to modify the {\it Clique-Coloring} tautology so that it becomes \emph{quasipolynomially} friendly in $Res$. This is done as follows:
\ben
\item consider $2k$-cliques vs. $k$ colorings;
\item add extension axioms for conjunctions up to $\log n$ that do not mix clique variables with coloring variables; this produces an {\bf NP}-disjoint pair equivalent to the canonical pair;
\item translate the quasipolynomial proofs of $PHP^{2k}_k$ in $Res(\log)$ into $Res$ proofs with the extension axioms.
\een
It does not matter that we only get quasipolynomially friendly disjuction, because the lower bound on the Clique-Coloring tautology is exponential.
\epri

The following concept is just a curiosity. 

\bdf
A disjunction $\phi(\vec p,\vec x)\vee\psi(\vec p,\vec y)$ is \emph{strongly friendly} for a proof system $P$ if
\ben
\item it is friendly and
\item for every $\vec a\in\{0,1\}^n$ such that both $\phi(\vec a,\vec x)$ and $\psi(\vec a,\vec y)$ are tautologies, both formulas are hard for $P$.
\een
\edf
A \emph{trivial strongly friendly disjunction} is one that has a polynomial size proof and one term is equivalent to the negation of the other.

\bt
There exists a nontrivial strongly friendly disjunction for Resolution.
\et

\bprf
We will use the fact that there are two polynomials $p_1$ and $p_2$ such that
\ben
\item if $\psi$ is a nontautlogy of length $n$ and $m\geq n$, then $\neg\pf_{Res,m,n}(\vec z,\lceil\psi\rceil)$ has a $Res$-proof of length $p_1(m)$, (polynomial upper bound on the local reflection principle for nontautologies, Proposition~\ref{p39});
\item if $\psi$ is a tautology and $m\geq p_2(n)$, then $\neg\pf_{Res,m,n}(\vec z,\lceil\psi\rceil)$ has only exponentially long proofs (Lemma~\ref{l39}). 
\een
Let $p$ be the maximum of $p_1$ and $p_2$. In the rest of the proof we will omit the subscript $Res$.

Our strongly friendly disjunction is: 
\[
\neg \pf_{m',n'}(\vec x,\lceil\neg \pf_{m,n}(\vec z,\vec y)\rceil(\vec y))\vee 
sat_{n'}(\lceil\neg \pf_{m,n}(\vec z,\vec y)\rceil(\vec y),\vec u),
\]
where $m=p(n)$, $n'=|\lceil\neg \pf_{m,n}(\vec z,\vec y)\rceil(\vec y)|$ and $m'=p(m)$. 
The expression $\lceil\neg \pf_{m,n}(\vec z,\vec y)\rceil(\vec y)$ has a similar meaning as $\lceil\neg Clique_{k+1}(\vec x,\vec y)\rceil(\vec y)$ in the proof of the lower bound on the reflection principle for $CP$. Namely, it is a string of constants and variables $\vec y$ such that if we substitute a code of a formula $\psi$ for $\vec y$, we get a code of the formula $\neg \pf_{m,n}(\vec z,\lceil\psi\rceil)$. The disjunction is a restriction of the reflection formula to a certain type of formulas, hence it is friendly. It remains to prove that if we instantiate $\vec y$ by plugging in some formula $\psi$ and obtaining
\[
\neg \pf_{m',n'}(\vec x,\lceil\neg \pf_{m,n}(\vec z,\lceil\psi\rceil)\rceil)\vee 
sat_{n'}(\lceil\neg \pf_{m,n}(\vec z,\lceil\psi\rceil)\rceil,\vec u),
\]
then either one of the two disjuncts 
is a nontautology, or both are hard for Resolution. We will consider three cases. Since $\psi$ is fixed, we can simplify the second term and get
\[
\neg \pf_{m',n'}(\vec x,\lceil\neg \pf_{m,n}(\vec z,\lceil\psi\rceil)\rceil)\vee 
\neg \pf_{m,n}(\vec z,\lceil\psi\rceil).
\]
\medskip\noindent
{\it Case (i), $\psi$ is a nontautology.} Then $\neg \pf_{m,n}(\vec z,\lceil\psi\rceil)$ has a proof of length $m'$, hence the first disjunct is a nontautology.

\medskip\noindent
{\it Case (ii), $\psi$ is a tautology and has a proof of length $\leq m$.} Then $\neg \pf_{m,n}(\vec z,\lceil\psi\rceil)$ is a nontautology.

\medskip\noindent
{\it Case (iii), $\psi$ is a tautology and has no proof of length $\leq m$.} Then
\bi
\item  $\neg \pf_{m,n}(\vec z,\lceil\psi\rceil)$ does not have a subexponential proof according to fact 2. above;
\item since $\neg \pf_{m,n}(\vec z,\lceil\psi\rceil)$ is a tautology, 
$\neg \pf_{m',n'}(\vec x,\lceil\neg \pf_{m,n}(\vec z,\lceil\psi\rceil)\rceil)$ also  does not have a subexponential proof according to fact 2.
\ei
\eprf

It would be interesting to find a nontrivial strongly friendly disjunction defined combinatorially. It could be one of the well known tautologies, such as the Clique-Coloring tautology, or the Broken Mosquito Screen tautology of~\cite{cook-haken}. If one disjunct defines an {\bf NP}-complete set, then such a pair would give us an alternative, syntax-free proof of non-automatibility of Resolution based on the condition {\bf P$\neq$NP}.

\subsubsection{Bounded depth Frege systems}

We will denote by $F_d$ the bounded depth Frege system formalized by the propositional sequent calculus with formulas restricted to depth $d$. Thus Resolution is~$F_0$. The following two statements are probably true:

\ben
\item \emph{$F_{d+1}$ provably $p$-proves the reflection principle for $F_d$,} 
\item \emph{$F_d$ provably $p$-proves the local reflection principle for $F_d$ for  non-tautologies.}
\een
The first appeared in Beckman et al.~\cite{BPT} with a proof idea. As for the second one, it seems that the proof that {$Res$ provably $p$-proves the reflection principle for non-tautologies for $Res$} can be generalized to all systems~$F_d$. The reason why believe that these are true fact is that the formulas expressing that a clause is satisfied by an assignment (general in the case of the reflection principle and specific in the case of the local reflection principle) have the appropriate depth, cf. the proof of Proposition~\ref{p39}. Another result in this vein is:
\ben
\item[3.] For every $k\geq 2$, $Res(k+1)$ provably $p$-proves the reflection principle for $Res(k)$, Atserias and Bonet~\cite{atserias-bonet}.
\een

It is also possible that other results about Resolution can be generalized to bounded depth Frege systems, but this will certainly require nontrivial work. We state two such generalizations as open problems.

\begin{problem}\label{p1}
Does $F_d$ $p$-prove its reflection principle for $d\geq 1$? The same question  for the local reflection principle for tautologies.
\end{problem}\begin{problem}\label{p2}
Is it true that for $d\geq 1$, every tautology $\phi$, and $m$ sufficiently large w.r.t. $n$,   $\neg\mbox{\it prf}_{F_d,m,n}(\vec x,\lceil\phi\rceil)$ 
doesn't have polynomial length $F_d$-proofs?
\end{problem}

The negative answer to Problem~\ref{p1} would show $F_d$ is weaker than $F_{d+1}$ on $DNF$ tautologies. A small superpolynomial separation has been proved by Impagliazzo and Kraj\'i\v cek~\cite{impagliazzo-krajicek},%
\footnote{In \cite{impagliazzo-krajicek} the result is only stated in terms of bounded arithmetical theories; cf.~\cite{krajicek-proof} Section 14.5 for the proof of the separation of proof systems.} 
but a  superquasipolynomial separation is still an open problem. A positive answer to Problem~\ref{p2} would imply that $F_d$ is non-automatable if {\bf P}$\neq${\bf NP} (but it would not refute 1).

To prove a superquasipolynomial separation of $F_d$ and $F_{d+1}$, it would suffice to prove a a superquasipolynomial lower bound on $F_d$-proofs of $\mbox{\it lrfn}_{P,\phi,m}$ for one sequence of tautologies $\phi_n$ that do not have quasipolynomial size $F_d$ proofs, e.g., $PHP_n$:
\[
\neg\mbox{\it prf}_{F_d,m,(n+1)n}(\vec x,\lceil PHP_n\rceil)\vee PHP_n,
\]
but we do not know how to prove lower bounds on the proof even if one only takes the first term of the disjunction. This is an instance of the fundamental problem: \emph{how difficult is it to prove a lower bound on the lengths of proofs?} In particular, we do not know the answer to the following:

\begin{problem}
Is there a $d\geq 1$ and a sequence of $DNF$ tautologies $\{\alpha_n\}$ such that a superpolynomial lower bound on the Resolution proofs of $\{\alpha_n\}$ can be proved in $F_d$ using polynomial size proofs? More precisely: is there a function $f$ growing more than polynomially such that tautologies $\neg\pf_{Res,f(n),n}(\vec x,\lceil\phi\rceil)$ have polynomial size $F_d$-proofs?
\end{problem}

\subsubsection{Frege and Circuit Frege systams}

We have already mentioned, Corollary~\ref{c3.2}, that all extensions of $CF$ provably $p$-prove their reflection principles. This result is just an easy generalization of Cook's proof of this fact for Extended Resolution~\cite{cook}. Buss constructed explicitly polynomial size proofs of the reflection principle for a Frege system in the Frege system~\cite{buss-reflection}. One can also prove this fact using the theory $\mbox{\it V}NC^1$ that is associated with Frege systems, see~\cite{cook-nguyen}.

The picture that emerges from what we know about Resolution, bounded depth Frege Systems on one side, and Frege, Circuit Frege systems and their extensions on the other is that these two kinds of proof systems have essentially different properties. However, so far we only know that the strong systems $p$-prove their reflection principles, while the weak ones probably do not. It would be interesting to find more differences. The most interesting question is whether the strong proof systems can prove lower bounds on their proofs. To prove that they cannot, of course, is hard, becasue we cannot prove any lower bounds on them, but it is conceivable that one can prove at least that such lower bounds are not provable in the theories associated with them, i.e., strengthen Proposition~\ref{p-t-unprov} below.


\section{Theories and proof systems}

We know that the fragments $G_i$ of the quantified propositional calculus prove their reflection principle using polynomial length proofs  (see~\cite{krajicek-pudlak-quantified}). We would like to argue that a proof system proves efficiently its reflection principle is rather a rule than an exception. To this end we will study theories associated with proof systems and proof systems associated with theories. 
Given an arithmetical theory $T$, we can associate two kinds of proof systems with~$T$. The first one, which we will call weak, may not always exist. The second one, which we will call strong, is always defined. Our terminology \emph{weak/strong proof system} is new. In  \cite{cook-nguyen}  the weak proof system of theory $T$ is called \emph{a proof system associated with $T$}; in \cite{krajicek-proof} it is called \emph{a proof system corresponding to $T$}. 
The strong proof system has been defined before (\cite{krajicek-pudlak-propositional} is, maybe, the first reference), but no name was given to it.

\subsection{The weak proof systems of theories}

\bdf\label{d-1}
Let $T$ be a theory and $P$ a proof system. We will say that
$P$ \emph{provably $p$-simulates} $T$ if for every $\forall\Sigma^b_0$ sentence $A$, $P$ provably $p$-proves propositions $[\![A]\!]_n$.
\edf



\bdf
We say that $P$ is a \emph{weak proof} system of an arithmetical theory $T$ if
\ben
\item $T$ proves the soundness of $P$, and  
\item $P$ provably $p$-simulates $T$.
\een
\edf
Given a theory we can always define a proof $P$ sytem that provably $p$-simulates $T$ and vice versa, given a proof system, we can define a theory that proves the reflection principle (soundness) of $P$, but these constructions in general do not ensure that both properties hold simultaneously.
The most important fact concerning weak systems of theories is the following theorem.

\bt[\cite{cook,krajicek-pudlak-quantified}]\label{t4.1}
Let $P$ be a weak proof system of a theory $T$. Then $P$ is the strongest proof system whose soundness is provable in $T$, i.e., every proof systems $P'$ whose soundness is provable in $T$ can be polynomially simulated by $P$. Moreover, the latter fact is provable in~$S^1_2$.
\et
\bprf
The following argument can be formalized in $S^1_2$. 
Let $P$ be a weak proof system of a theory $T$. If $T$ proves the soundness of $P'$, which is the sentence $\Rfn_{P'}$, then $P$ proves the propositional instances of $\Rfn_{P'}$, which are propositional formulation of the reflection principle $\rfn_{P',m,n}$. Since we are arguing in $S^1_2$, the $P$-proofs of $\rfn_{P',m,n}$ can be constructed in polynomial time in $n$ and $m$. If we want to prove a proposition $\phi$ in $P$ and we are given a $P'$-proof $D$ of $\phi$, we only need to substitute $\phi$ and $D$ into $\rfn_{P',m,n}$ where $m=|\phi|,n=|D|$.
\eprf
\begin{corollary}
If $P$ and $P'$ are weak proof systems of a theory $T$, then they are polynomially equivalent.
\end{corollary}
According to this corollary, a weak proof system for a theory $T$ is determined up to polynomial simulation. We will denote by $P_T$ one of these weak proof systems (when they exist). Note also that Theorem~\ref{t4.1} gives us an equivalent definition of the weak propositional proof system of a theory.

Let $A\in\forall\Sigma^b_0$. We will denote by $CF+\{[\![A]\!]_n\}$ the Circuit Frege proof system extended with axiom schemas $\{[\![A]\!]_n\}$, which means that the system can use any proposition of the form $[\![A]\!]_n(x_1/\beta_1\dts x_m/\beta_m)$ as an axiom, where $x_1\dts x_m$ are the propositional variables of $[\![A]\!]_n$ and $x_i/\beta_i$ denotes the substitution of a proposition $\beta_i$ for variable $x_i$.

\bt[\cite{krajicek-pudlak-quantified}]\label{t43}
~
\ben
\item For every true $\forall\Sigma^b_0$ sentence $A$ sentence, $CF+\{[\![A]\!]_n\}$ is a weak proof system of $S^1_2+A$.
\item A proof system $P$ that extends $CF$ is a weak proof system of $S^1_2+\mbox{\it Rfn}_P$ iff $P$ provably $p$-proves its reflection principles.
\een
\et
\bprf
1. First we show that $S^1_2+A$ proves the soundness of  $CF+\{[\![A]\!]_n\}$. It proves the soundness of the (substitution instances of) axioms $\{[\![A]\!]_n\}$ using $A$. Then the argument is the same as for $CF$ alone. Since $S^1_2+A$ proves that the rules preserve soundness, $PIND$-$\Pi^b_1$ implies that the proof system is sound.

Now we prove that $CF+\{[\![A]\!]_n\}$ simulates $S^1_2+A$. Let $A:=\forall y.\phi(y)$, where $\phi\in\Sigma^b_0$. Suppose $S^1_2+A$ proves $\forall x.\psi(x)$ for some $\psi\in\Sigma^b_0$. By Buss's theorem,
\[
S^1_2\vdash\ \forall x(\phi(t(x))\to \psi(x)),
\]
where $t(x)$ represents a polynomial time computable function. Since $CF$ simulates $S^1_2$, we know that $S^1_2$ proves that $[\![\phi(f(x)\to\psi(x)]\!]_n$ are provable in $CF$. Since the substitution instances of $[\![A]\!]_n$ are axioms of $CF+\{[\![A]\!]_n\}$, $S^1_2$ proves that $[\![\phi(f(x))]\!]_n$ are provable in $CF+\{[\![A]\!]_n\}$. Hence it also proves that $\{[\![\psi(x)]\!]_n\}$ are provable in $CF+\{[\![A]\!]_n\}$.

\medskip
2. ($\Rightarrow$) If $P$ is a weak proof system of $S^1_2+\mbox{\it Rfn}_P$, then in particular $S^1_2$ proves that $P$ proves $[\![\mbox{\it Rfn}_P]\!]_n$ for all $n$. These propositions are equivalent to the propositions $\rfn_{P,k,l}(\vec x,\vec y,\vec z)$ expressing the reflection principle in the propositional calculus. 

($\Leftarrow$) Suppose that $S^1_2$ proves that $P$ proves its reflection principle. Trivially, $S^1_2+\mbox{\it Rfn}_P$ proves the soundness of $P$. It remains to show that $P$ simulates $S^1_2+\mbox{\it Rfn}_P$. 
The argument is very similar to the one in the proof in part 1 of the simulation of $S^1_2+A$ by $CF+\{[\![A]\!]_n\}$. The difference is that we now use $\mbox{\it Rfn}_P$ instead of $A$, and we use the fact that $P$ proves the propositional instances of the reflection principle instead of $\{[\![A]\!]_n\}$.
\eprf

\begin{corollary}\label{c4.4}
~
\ben
\item For every $\forall\Sigma^b_0$ sentence $A$, $S^1_2+A$ has a weak proof system.
\item For every proof system $P$ that extends $CF$, $P+\{\rfn_{P,m,n}\}_{m,n}$ is a weak proof system of some theory. Furthermore, $P+\{\rfn_{P,m,n}\}_{m,n}$ provably $p$-proves propositions expressing its reflection principle.
\een
\end{corollary}
\bprf
1. Follows immediately from Theorem~\ref{t43},~1.

2. By Theorem~\ref{t43},~1, $CF+\{\rfn_{P,m,n}\}_{m,n}$ is a weak proof system of $S^1_2+\Rfn_P$. Further, $CF+\{\rfn_{P,m,n}\}_{m,n}$ polynomially simulates $P+\{\rfn_{P,m,n}\}_{m,n}$, because $CF+\{\rfn_{P,m,n}\}_{m,n}$ polynomially simulates $P$; the reverse simulation is trivial. Thus $P+\{\rfn_{P,m,n}\}_{m,n}$ is a weak system of $S^1_2+\Rfn_P$. By Theorem~\ref{t43},~2,  $P+\{\rfn_{P,m,n}\}_{m,n}$ provably $p$-proves its reflection principle.
\eprf

Now we are almost ready to prove that nonlinear lower bounds on $P$-proofs as stated in Proposition~\ref{p-unprov} are not provably $p$-provable in $P$. It only remains to prove the unprovability in the theory associated with $P$.

\bprl{p-t-unprov}
Let $P$ be a weak proof system of a theory $T\supseteq S^1_2$, let $Q$ be a proof system which extends $CF$,
and let $\{\alpha_n\}$ be a polynomial time computable sequence of tautologies with $|\alpha_n|=n$. Then there exists a constant $c$ such that for every $n_0$, 
if $T$ proves a lower bound $>cn$ for $n\geq n_0$ on $Q$-proofs of $\alpha_n$, then $T$ proves the soundness of $Q$, hence $P$ polynomially simulates~$Q$.%
\footnote{Note that this does not exclude the possibility that $T$ proves a superlinear lower bound on $Q$ that is stronger than $P$ \emph{assuming in $T$ that $Q$ is sound.}}
\epr
\bprf
Lemma~\ref{l5.1} can be formalized in $S^1_2$. Thus there exists a constant $c_1$ depending only on $Q$ such that $T$ proves:
\ben
\item[(*)]{\it if $Q$ is not sound, then there exists an $r$ such that all propositions have $Q$-proofs of length $\leq c_1n+r$. }
\een
Let $c>c_1$ and suppose $T$ proves a lower bound $>cn$ for $n\geq n_0$ on $Q$-proofs of $\alpha_n$. Then $T$ proves that for every $r$ there exists an $n$ such that $cn\geq c_1n+r$. 
Hence according to (*), $T$ proves that $Q$ is sound.
\eprf

\bprf[of Proposition \ref{p-unprov}]
Let proof systems $P$, $Q$ and tautologies $\{\alpha_n\}$ be given, and suppose 
$P$ provably $p$-proves its reflection principle and a lower bound $>cn$ for $n\geq n_0$ on $Q$-proofs of $\alpha_n$. Let $T$ be $S^1_2+\Rfn_P$. By Theorem~\ref{t43}, $P$ is a weak proof system of $T$. Our assumption is that $S^1_2$ proves that $P$ proves  a lower bound $>cn$ for $n\geq n_0$ on $Q$-proofs of $\{\alpha_n\}$. Since $T$ proves the soundness of $P$, it also proves a lower bound $>cn$ for $n\geq n_0$ on $Q$-proofs of $\{\alpha_n\}$. Hence, by the previous proposition, $P$ polynomially simulates $Q$.
\eprf

One can prove unprovability of superpolynomial lower bounds on $EF$ in $PV_1$, which is a theory slightly weaker than $S^1_2$ and for which $EF$ is the weak proof system. 

\bt[\cite{krajicek-pudlak-models}, \cite{krajicek-proof} Section 20.1]
For every function $f:\N\to\N$ such that for every $c\in\N$, $PV_1$ proves that $f(n)$ eventually dominates $n^c$, $PV_1$ does not prove the sentence
\[
\forall x\exists y (|x|<|y|\wedge Taut(y)\wedge\forall z(|z|\leq f(|y|)\to\neg\Prf_{EF}(z,y))).
\]
\et
This can certainly be generalized to stronger theories and proof systems, but we do not see a way to derive from it a statement about unprovability in propositional proof systems.

\subsection{The strong proof systems of theories}

\bdf
Let $T$ be a consistent theory extending Robinson's Arithmetic with the set of axioms decidable in polynomial time. The \emph{strong proof system of} $T$, denoted by $Q_T$, is the propositional proof system where $d$ is a proof of a proposition $\phi$, if $d$ is a $T$-proof of the sentence $Taut(\lceil\phi\rceil)$.
\edf
Since Robinson's Arithmetic proves all true bounded sentences and $T$ is consistent, $Q_T$ is complete and sound. There is an \emph{alternative definition} that works for all theories that have infinite models. This is based on a more direct translation of propositional formulas into first-order sentences. Let $\phi(p_1\dts p_n)$ be a propositional formula where $p_1\dts p_n$ are the propositional variables of $\phi$. We define its first-order translation $tr_\phi$ by taking $n+1$ distinct first-order variables $x_1\dts x_n,y$ and putting
\[
tr_\phi:= \forall x_1\dots x_ny.\phi(x_1=y\dts x_n=y).
\]

{\it Exercise.} Construct polynomial size proofs of the sentences $tr_\phi\equiv Taut(\lceil\phi\rceil)$ in Robinson's Arithmetic.

\begin{fact}
Let $T$ be theory such that $P_T$ is defined. Then $Q_T$ polynomially simulates $P_T$.
\end{fact}
\bprf
As $T$ proves $\Rfn_P$, given a proposition $\phi$ and its $P$-proof, we can construct in polynomial time a $T$-proof of $Taut(\lceil\phi\rceil)$.
\eprf

One reason for calling $Q_T$ strong is the following fact.

\begin{fact}\label{f3}
$T$ does not prove the soundness of $Q_T$.
\end{fact}
\bprf
By G\"odel's 2nd incompleteness theorem.
\eprf
This fact, however, does not exclude the possibility that $Q_T$ is not strong, e.g., that it is polynomially equivalent to $P_T$, although we consider this possibility unlikely. On the contrary, it seems reasonable to conjecture that for every theory $T$ that possesses a weak proof system, $Q_T$ is strictly stronger than $P_T$. The following corollary was suggested by J.~Pich.
\begin{corollary}
$T$ does not prove superlinear lower bounds on $Q_T$-proofs of any sequence of tautologies~$\{\alpha_n\}$.
\end{corollary}
\bprf
From Fact~\ref{f3} and Proposition~\ref{p-t-unprov}.
\eprf

\subsection{Consistency statements}

For a theory $T$ whose set of axioms is polynomial time decidable, we denote by $Con_T(x)$ an arithmetical sentence that formalizes the statement that there is no proof of contradiction in $T$ whose length is $\leq x$. In this formula, $x$ is a variable ranging over natural numbers. Thus the consistency of $T$ can be expressed by $\forall x.Con_T(x)$.
For a fixed $n\in\N$, we denote by $\bar n$ a closed arithmetical term of length $O(\log n)$ whose value is~$n$. Hence $Con_T(\bar n)$ is a sentence of length $O(\log n)$ formalizing the statement that there is no proof of contradiction of length~$\leq n$. 

\bt[\cite{finitistic,improved}]\label{t4.5}
Let $T$ be a finitely axiomatized sequential theory. Then there exists a polynomial time computable sequence of $T$-proofs $\{D_n\}$ such that $D_n$ is a proof of $Con_T(\bar n)$ for $n=1,2,\dots$. Moreover, it is provable in $S^1_2$ that for all $n$, $D_n$ is a  $T$-proof of  $Con_T(\bar n)$.
\et
In \cite{finitistic,improved} we constructed a sequence of $T$-proofs $\{D_n\}$ of $Con_T(\bar n)$ of polynomial length in~$n$. One can easily check that the construction can be done in polynomial time and this is also provable in~$S^1_2$.
In contrast to Theorem~\ref{t4.5}, if $\{E_n\}$ is a sequence of $T$-proofs of $Con_{T+Con_T}(\bar n)$, then this fact is not provable even in $T+Con_T$.

\bpr
Let $T$ be a consistent computably axiomatized theory containing $S^1_2$. Then for no sequence $\{E_n\}$, $T+Con_T$ proves that for all $n$,  $E_n$ is a  $T$-proof of $Con_{T+Con_T}(\bar n)$.
\epr
\bprf
By way of contradiction, suppose $\{E_n\}$ is such a sequence. Since $T+Con_T$ proves the consistency of $T$, it also proves the uniform $\Pi_1$ reflection principle for $T$. Since it proves that $\{E_n\}$ is a sequence of $T$-proofs of $Con_{T+Con_T}(\bar n)$, it also proves, using the reflection principle, $\forall x.Con_{T+Con_T}(x)$, which is equivalent to $Con_{T+Con_T}$. This is in contradiction with the Second Incompleteness Theorem.
\eprf

If $T$ contains $S^1_2$ and $T+Con_T$ is consistent, then there exists a computable sequence $\{E_n\}$ of $T$-proofs of $Con_{T+Con_t}(\bar n)$. It is a formalization of the brute-force search for contradiction and the proofs have exponential lengths in $n$. We have conjectured that there are no such proofs of polynomial length, see~\cite{finitistic}. A slightly weaker conjecture says that such proofs cannot be constructed in polynomial time. 
In the next section we will show a link with the conjecture that the strong system of a theory $T$ is strictly stronger than the weak one.

\subsection{Strong vs. weak proof systems}

\bt[essentially, Kraj\'i\v cek's Lemma 12.8.2 of \cite{krajicek-proof}]\label{t3.1}
Let $T=S^1_2+A$, where $A$ is a true $\Pi_1$ sentence. Then $Q_T$ is a weak proof system of $T+Con_T$, i.e.,
\[
Q_T\equiv P_{T+Con_T}.
\]
\et
\bprf
Clearly, $T+Con_T$ proves the consistency of $T$, hence also the soundness of $Q_T$. It remains to prove that provably in $S^1_2$, $Q_T$ simulates $T+Con_T$. Let a $\Sigma^b_0$ formula $\theta(x)$ be given and suppose that $T+Con_T$ proves $\forall x.\theta(x)$. Then
\[
T\vdash(\forall y.Con_T(|y|))\to\forall x.\theta(x).
\]
By Parikh's theorem, there exists a polynomial $p$ such that 
\[
T\vdash \forall x\exists y(|y|\leq p(|x|)\wedge (Con_T(|y|)\to\forall x.\theta(x))).
\]
This implies
\[
T\vdash \forall x(Con_T(p(|x|))\to\theta(x)).
\]
Since $T\vdash\ u\leq v\to (Con_T(v)\to Con_T(u))$, we have 
\[
T\vdash \forall x(Con_T(p(|x|))\to\forall z(|z|\leq|x|\to\theta(z)).
\]
Let $D$ be a $T$-proof of this sentence. By Theorem~\ref{t4.5}, provably in $S^1_2$, there exists an algorithm that in time polynomial in $n$ constructs a $T$-proof $D_n$ of $Con_T(n)$.  Thus
$S^1_2$  proves that for all $n$,  $D+D_{p(n)}$ is a $T$-proof of $\forall z(|z|\leq \bar n\to\theta(z))$. We also have 
\[
T\vdash\forall z(|z|\leq \bar n\to\theta(z))\equiv Taut([\![\theta]\!]_n)
\]
with a polynomial size proof and provably in $S^1_2$. 
Thus $Q_T$ provably $p$-simulates $T+Con_T$. 
\eprf

\begin{corollary}
Let $T=S^1_2+A$, where $A$ is a $\forall\Sigma^b_0$ true sentence, and let $S=T+Con_T$. Then TFAE:
\ben
\item $T$-proofs of $Con_S(\bar n)$ can be constructed in time polynomial in $n$.
\item The weak system of $S$ polynomially simulates the strong system of $S$.
\een
\end{corollary}

\bprf
We will use the observation that $Con_S(\bar n)$ is equivalent to the statement that $con_{Q_S,n}(\vec x)$ is a tautology, i.e.,
\bel{e-con}
S^1_2\vdash\ Con_S(\bar n)\ \equiv\ Taut(con_{Q_S,n}).
\ee

\medskip\noindent
1$\Rightarrow$2. 
Suppose $T$-proofs of $Con_S(\bar n)$ can be constructed in polynomial time. This is equivalent to constructing $Q_T$-proofs of $con_{Q_S,n}$ in polynomial time. By Corollary~\ref{c3.2} we get $Q_T$-proofs of the reflection principle for $Q_S$ in polynomial time. Hence, given a $Q_S$-proof of some $\phi$, we get a $Q_T$ proof of $\phi$ in polynomial time. Thus $Q_T$ polynomially simulates $Q_S$. By Theorem~\ref{t3.1}, this implies that $P_S$ polynomially simulates $Q_S$.

\medskip\noindent
2$\Rightarrow$1. By Theorem~\ref{t4.5}, there exists an algorithm that in time polynomial in $n$ constructs $S$-proofs $Con_S(n)$. Hence by~(\ref{e-con}) above, one can construct in polynomial time $Q_S$-proofs of $Taut(con_{Q_S,n})$. Assuming $P_S$ polynomially simulates $Q_S$, we get $P_S$-proofs of $Taut(con_{Q_S,n})$. By Theorem~\ref{t3.1}, this implies that we get $Q_T$-proofs of $Taut(con_{Q_S,n})$ in polynomial time. Using~(\ref{e-con}) once again, we get $T$-proofs of $Con_S(\bar n)$.
\eprf

One can easily check that a version of this corollary with polynomial lengths of proofs instead of polynomial time algorithms is also true.
It seems that in general $Q_T$ is \emph{much} stronger than $P_T$. Let's consider an example.

\medskip {\it Example.} Let $T$ be $S^1_2$. Then, by~\cite{cook,buss}, $P_T$ is the Extended Frege proof system. It is well-known that $S^1_2$ interprets the entire bounded arithmetic $T_2$ on an initial segment of the natural numbers. $T_2$ proves the soundness of all fragments $G_i$ of the quantified propositional sequent calculus~\cite{krajicek-pudlak-quantified}. From this, one can easily deduce that the strong system of $S^1_2$ polynomially simulates all $G_i$s. But in fact, the strong system of $S^1_2$ simulates (apparently) much stronger proof systems. Kraj\'i\v cek defined a construction that, from a proof system $P$, produces the \emph{implicit} $P$, denoted by $iP$, and this can be iterated (see~\cite{krajicek-implicit}). The implicitation construction seems to always produce a stronger system. E.g., $iEF$ polynomially simulates $G$, the quantified propositional sequent calculus, while we believe that $EF$ does not. The strong proof system of $S^1_2$ polynomially simulates all iterated implicitations $i_kEF$ of $EF$. This is because
\ben
\item $I\Delta_0+Exp$ proves the soundness of all $i_kEF$, and
\item if  $I\Delta_0+Exp$ proves a $\forall\Sigma^b_0$ sentence $\phi$, then there exists a formula $\alpha(x)$ with one free variable that defines an initial segment of numbers closed under the successor, $+$, and $\times$ in $S^1_2$ such that $S^1_2$ proves the relativization of $\phi$ to $\alpha$; 
see~\cite{hajek-pudlak}. 
\een

\begin{problem}
Characterize the strong proof system of $S^1_2$.
\end{problem}

We conjecture that $i_\infty EF$, a suitably formalized union of all $i_kEF$s, is polynomially equivalent to the strong proof system of $S^1_2$. We do have some ideas how to prove this conjecture, but for the time being it is an open problem.



\end{document}